\documentclass[11pt]{article}
\usepackage{latexsym}
 \usepackage{amssymb}
 \usepackage{amsmath}
 \usepackage{amsthm}
 \usepackage{graphicx}

 \newcommand{ \un }{\mathbb{I} }
 \newcommand{ \p }{\mathbb{P} }
 \newcommand{ \pa }{\mathbb{P}^{\alpha} }
  \newcommand{ \pao }{\mathbb{P}^{\alpha}_{0} }
\newcommand{ \pam }{\mathbb{P}^{\alpha}_{m_n}}
 \newcommand{ \E }{\mathbb{E}}
 \newcommand{ \Ea }{\mathbb{E}^{\alpha}}
 \newcommand{ \Eam }{\mathbb{E}^{\alpha}_{m_n}}

\newcommand{ \Ln }{\mathbb{L}_{n} }

\newcommand{ \F }{ \mathbb{F} }

 \newcommand{ \Z }{ \mathbb{Z} }
 \newcommand{\N}{ \mathbb{N} }
 
 \newcommand{ \V }{\textrm{Var} }

 \newcommand{ \Ct }{ \mathcal{C} }
 
 \newcommand{ \f }{ \mathcal{F} }

 \newcommand{ \tm }{ m }
 \newcommand{ \tM }{ M }
\newcommand{ \tmo }{ m_n }

\newcommand{\ra}{R^{\alpha}}

\newcommand{ \A}{ \mathcal{A} }

 \newcommand{ \lo }{ \mathcal{L} }

\newtheorem{The}{{\bf Theorem}}[section]
\theoremstyle{definition}
\newtheorem{Def}[The]{{\bf Definition}}
\theoremstyle{plain}
 \newtheorem{Lem}[The]{Lemma}
 \newtheorem{Cor}[The]{\bf Corollary}
 \newtheorem{Pro}[The]{\bf Proposition}
 \theoremstyle{definition}
\newtheorem{Rem}[The]{{\bf Remark}}

\newcommand{ \sa}[2]{E^{\alpha}_{#1}(#2)}

 \newenvironment{Pre}{\noindent \textbf{Proof.} \\ }{$\
 \blacksquare$}

\makeatletter

\@addtoreset{equation}{section} \makeatother

\title{Sinai's walk : a statistical aspect 
\\ \vspace{1cm}
 \large{Pierre Andreoletti} \footnote{Laboratoire MAPMO - C.N.R.S. UMR 6628 - F\'ed\'eration Denis-Poisson, Universit\'e d'Orl\'eans, 
(Orl\'eans France). \newline \vspace{0.1cm}  $\quad$  MSC 2000 60G50; 60J55; 62G05. \newline \vspace{0.5cm} \textit{Key words and phrases :  random environment, random walk,
Sinai's regime, non-parametric statistics} } \textrm{ }    }

\begin{document}

\bibliographystyle{unsrtnat}
\maketitle

%\noindent  $^*$ Laboratoire MAPMO - C.N.R.S. UMR 6628 - F\'ed\'eration Denis-Poisson, Universit\'e d'Orl\'eans.

\noindent \\ \textbf{Abstract:} We consider Sinai's random walk in
random environment. We prove that the logarithm of the local time is a good estimator of the random potential associated to the random environment. We give a constructive method allowing us to built the random environment from a single trajectory of the random walk.

%Thanks to this result we show, with the help of numerical simulations, that we can reconstruct  a significant part of the random potential and show an interval of confidence.

%\bibliographystyle{unsrt}

\section{Introduction and results}

In this paper we are interested in Sinai's walk i.e., a one
dimensional random walk in random environment with three
conditions on the random environment: two necessaries hypothesis
to get a recurrent process (see \cite{Solomon}) which is not a
simple random walk and an hypothesis of regularity which allows us
to have a good control on the fluctuations of the random
environment.

The asymptotic behavior of such walk has been understood by
\cite{Sinai} : this walk is sub-diffusive and at an instant $n$ it is localized in the neighborhood of a well
defined point of the lattice.  It is well known, see (Zeitouni [2001] for a survey) that this behavior is strongly dependent of the random environment or, equivalently, by the associated random potential defined Section 1.2.

The question we solve here is the following: given a single trajectory of a random walk  $(X_l, 1 \leq lÊ\leq n)$ where the time $n$ is fixed, can we estimate the trajectory of the random potential where the walk lives ? Let us remark that the law of this potential is unknown as-well. \\
\noindent In their paper, \cite{AdeEnr} are interested in the question of the distribution of the random environment that could be deduced from a single trajectory of the walk, on the other hand, our purpose is to get an approximation of the trajectory of  the random potential.  \\
In the paper \cite{Monasson} the authors are interested in a method to predict the sequence of DNA molecules. 
They model the unzipping of the molecule as a one-dimensional biased random walk for the fork position (number of open base pair) $k$ in this landscape. The elementary opening $(k \rightarrow k+1)$ and closing $(k \rightarrow k-1)$ transitions happen with a probability that depends on the unknown sequence. This probability of transition follows an Arrh\'enius law wich is closed to the one we discuss here. The question they answer is: given an unzipping signal can we predict the uniziping sequence ? 
%Their method is based on  a one dimensional random walk in an elementary unknown environment. 
Their approach is based on a Bayesian inference method which gives very good probabilities of prediction for a large amount of data. This  means, in term of the walk, several trajectory on the same environment. \\
Our approach is purely probailistic, it is based on good properties of the \textit{local time} of the random walk which is the amount of time the walk spends on the points of the lattice. We treat a general case with a very few information on the random environment. 
We are able to reconstruct the random potential in a significant interval where the walk spends most of its time. 
%The key point of this paper is that if we impose to the local time to be large enough but negligible comparing to the maximum of the local time then this will directly implies conditions on the random potential. Notice that, this aspect of looking the random environment from the point of view of the walk have already been studied for a different purpose in for example \cite{Kesten1}.  
Our proof is based on the results of \cite{Pierre1}, in particular in a weak law of large number for the local time on the point of localization of the walk.
%In fact it is well known that we can reconstruct the random potential using the local time. 

The largest part of this paper is devoted to the proof of a theoretical result (Theorem \ref{th1}), we also  present,  at the end of the document, numerical simulations to illustrate our result. We give the main steps of the  algorithm we use to rebuilt the random potential only by considering a trajectory of the walk. As an introduction we would like to comment one of these simulations:
\begin{figure}[h]
\begin{center}
\includegraphics[width=10cm,height=5cm]{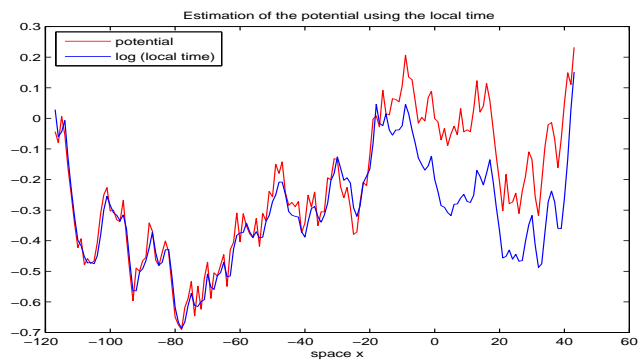}
\caption{The logarithm of the local time (in blue) and the random potential (in red)}
\end{center}
\end{figure}

In blue we have represented the logarithm of the local time and in red the potential associated to the random environment. 
First, remark that we get a good approximation on a large neighborhood of the bottom of the valley around the coordinate  -80. Outside this neighborhood and especially after the coordinate -20, the approximation is not precise at all. We will explain this phenomena by the fact that after the walk has reached the bottom of the valley, the walk will not return frequently to the points with coordinate larger than -20, so we lose information for this part of the latice.

Our method of estimation give us two crucial information: a confidence interval for the differencies of potential in sup-norm, on an observable set of sites ``sufficiently'' visited by the walk and a localization result for the bottom of the valley linked with the hitting time of the maximum of the local times.
First we need to define the process:

\subsection{Definition of Sinai's walk}

Let $\alpha =(\alpha_i,i\in \Z)$ be a sequence of i.i.d. random
variables taking values in $(0,1)$ defined on the probability
space $(\Omega_1,\f_1,Q)$, this sequence will be called random
environment. A random walk in random environment (denoted
R.W.R.E.) $(X_n,n
\in\N)$ is a sequence of random variable taking value in $\Z$, defined on $( \Omega,\f,\p)$ such that \\
$ \bullet $ for every fixed environment $\alpha$, $(X_n,n\in \N)$
 is a Markov chain with the following transition probabilities, for
 all $n\geq 1$ and $i\in \Z$
 \begin{eqnarray}
& & \p^{\alpha}\left[X_n=i+1|X_{n-1}=i\right]=\alpha_i, \label{mt} \\
& & \p^{\alpha}\left[X_n=i-1|X_{n-1}=i\right]=1-\alpha_i \equiv
\beta_i. \nonumber
\end{eqnarray}
We denote $(\Omega_2,\f_2,\pa)$ the probability space
associated to this Markov chain. \\
 $\bullet$ $\Omega = \Omega_1 \times \Omega_2$, $\forall A_1 \in \f_1$ and $\forall A_2 \in \f_2$,
$\p\left[A_1\times
A_2\right]=\int_{A_1}Q(dw_1)\int_{A_2}\p^{\alpha(w_1)}(dw_2)$.

\noindent \\ The probability measure $\pa\left[\left.
.\right|X_0=a \right]$  will be  denoted $\pa_a\left[.\right]$,
 the expectation associated to $\pa_a$: $\Ea_a$, and the expectation associated to $Q$:
 $\E_Q$.

\noindent \\ Now we introduce the hypothesis we will use in all
this work. The two following hypothesis are the necessaries
hypothesis
\begin{eqnarray}
 \E_Q\left[ \log
\frac{1-\alpha_0}{\alpha_0}\right]=0 , \label{hyp1bb} \label{hyp1}
\end{eqnarray}
\begin{eqnarray}
\V_Q\left[ \log \frac{1-\alpha_0}{\alpha_0}\right]\equiv \sigma^2
>0 . \label{hyp0}
\end{eqnarray}
 \cite{Solomon} shows that under \ref{hyp1} the
process $(X_n,n\in \N)$ is $\p$ almost surely recurrent and
\ref{hyp0} implies that the model is not reduced to the simple
random walk. In addition to \ref{hyp1} and \ref{hyp0} we will
consider the following hypothesis of regularity, there exists $0<
\eta_0 < 1/2$ such that
\begin{eqnarray}
& & \sup \left\{x,\ Q\left[\alpha_0 \geq x \right]=1\right\}= \sup
\left\{x,\ Q\left[\alpha_0 \leq 1-x \right]=1\right\} \geq \eta_0.
\label{hyp4}
\end{eqnarray}

\noindent We call \textit{Sinai's random walk} the random walk in
random environment previously defined with the three hypothesis
\ref{hyp1}, \ref{hyp0} and \ref{hyp4}.

\noindent
Let us define the local time $\lo$, at $k$ $(k\in \Z)$ within the
interval of time $[1,T]$ ($T \in \N^*$) of $(X_n,n\in \N)$
\begin{eqnarray}
\lo\left(k,T\right) \equiv \sum_{i=1}^T \un_{\{X_i=k\}} .
\end{eqnarray}
$\un$ is the indicator function ($k$ and $T$ can be deterministic or random variables). Let $V\subset \Z$, we denote
\begin{eqnarray}
\lo\left(V,T\right) \equiv \sum_{j \in V} \lo\left(j,T\right)
=\sum_{i=1}^T\sum_{j \in V} \un_{\{X_i=j\}} .
\end{eqnarray}
 
\noindent To end, we define the following random variables 

\begin{eqnarray}
& & \lo^*(n)=\max_{k\in \Z}\left(\lo(k,n)\right) ,\
\F_{n}=\left\{k \in \Z,\ \lo(k,n)=\lo^*(n) \right\} , \label{ref1.7} \\ 
& & k^*=\inf\{|k|, k\in \F_{n} \}
%&& Y_n=\inf_{x \in \Z}\min\left\{k>0\ :\ \lo([x-k,x+k],n) \geq n/2\right\}.
\end{eqnarray}

\noindent
 $\lo^*(n)$ is the maximum of the local times (for a given instant
$n$), $\F_{n}$ is the set of all the favourite sites and $k^{*}$ the smallest favorite site.
%$Y_{n}$ is the size of the interval where the walk spends more than a half of its %time.

\subsection{The random potential and the valleys}

From the random environment we define what we will call random potential, 
\noindent \\ Let
\begin{eqnarray}
\epsilon_i \equiv \log \frac{1-\alpha_i}{\alpha_i},\ i\in \Z,
\end{eqnarray}
define :
\begin{Def} \label{defpot2} The random potential $(S_m,\  m \in
\Z)$ associated to the random environment $\alpha$ is defined in the following way: 
\begin{eqnarray} 
 && S_{k}=\left\{ \begin{array}{ll} \sum_{1\leq i \leq k} \epsilon_i, & \textrm{ if }\ k >0, \\
  -\sum_{k+1 \leq i \leq 0} \epsilon_i , &   \textrm{ if }\  k<0  , \end{array} \right.  \nonumber \\
&& S_{0}=0. \nonumber
 \end{eqnarray}
\end{Def}
% \begin{eqnarray} \left\{ \begin{array}{ll} 
% && S_{k}-S_{j}= \sum_{j+1\leq i \leq k} \epsilon_i, \\
% && S_{0}=0, \nonumber \end{array} \right.
% \end{eqnarray}
%and symmetrically if $k<j$.

\begin{figure}[h]
\begin{center}
%\input{thfig4.pstex_t} \caption{} \label{thfignew1}
%\end{center}
%\end{figure}
%\begin{figure}[h]
%\begin{center}
\input{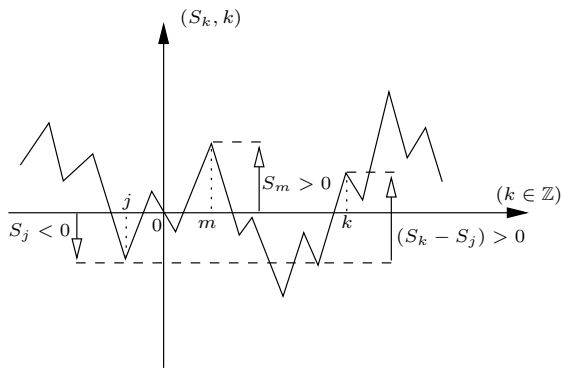} 
\caption{Trajectory of the random potential} \label{fig5}
\end{center}
\end{figure}

\begin{Def} \label{c2s2d1}
 We will say that the triplet $\{M',m,M''\}$ is a \textit{valley} if
 \begin{eqnarray}
& &  S_{M'}=\max_{M' \leq t \leq m} S_t ,  \\
& &  S_{M''}=\max_{m \leq t \leq
 \tilde{M''}}S_t ,\\
& & S_{m}=\min_{M' \leq t \leq M''}S_t \ \label{2eq58}.
 \end{eqnarray}
If $m$ is not unique we choose the one with the smallest absolute
value.
 \end{Def}

\begin{Def} \label{deprofvalb}
 We will call \textit{depth of the valley} $\{\tM',\tm,\tM''\}$ and we
 will denote it $d([M',M''])$ the quantity
\begin{eqnarray}
 \min(S_{M'}-S_{m},S_{M''}-S_{m})
 .
 \end{eqnarray}
 \end{Def}

\noindent  Now we define the operation of \textit{refinement}
 \begin{Def}
Let  $\{\tM',\tm,\tM''\}$ be a valley and let
  $\tM_1$ and $\tm_1$ be such that $\tm \leq \tM_1< \tm_1 \leq \tM''$
  and
 \begin{eqnarray}
 S_{\tM_1}-S_{\tm_1}=\max_{\tm \leq t' \leq t'' \leq
 \tM''}(S_{t'}-S_{t''}) .
 \end{eqnarray}
 We say that the couple $(\tm_1,\tM_1)$ is obtained by a \textit{right refinement} of $\{\tM',\tm,\tM''\}$. If the couple $(\tm_1,\tM_1)$ is not
 unique, we will take  the one such that $\tm_1$ and $\tM_1$ have the smallest  absolute value. In a similar way we
define the \textit{left refinement} operation.
 \end{Def}

\begin{figure}[h]
\begin{center}
%
%\input{thfig4.pstex_t} \caption{} \label{thfignew1}
%\end{center}
%\end{figure}
%\begin{figure}[h]
%\begin{center}
\input{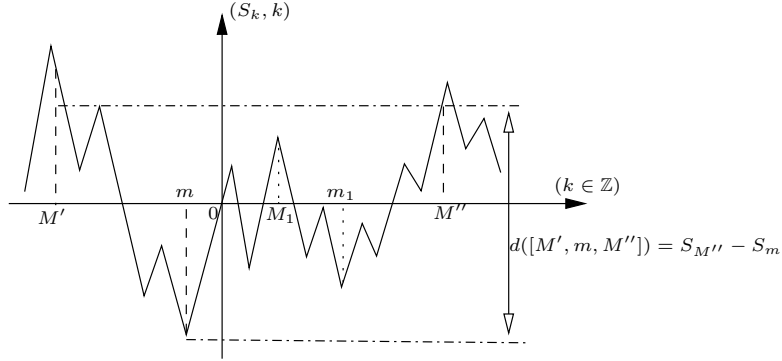}
 \caption{Depth of a valley and refinement operation} \label{fig4}
\end{center}
\end{figure}

\noindent \\  We denote $\log_2=\log \log $, in all this section we will suppose that $n$ is
large enough such that $\log_2 n$ is positive.

\begin{Def} \label{thdefval1b} Let $n>3$, $\gamma>0$, and $\Gamma_n \equiv \log n+ \gamma \log_2 n $, we say that a valley $\{\tM',\tm,\tM''\}$ contains $0$ and
is of depth larger than $\Gamma_n$ if and only if
\begin{enumerate}
\item  $ 0 \in [\tM',\tM'']$, \item $d\left([\tM',\tM'']\right)
\geq \Gamma_n $ , \item if $\tm<0,\ S_{\tM''}-\max_{ \tm \leq t
\leq 0
}\left(S_t\right) \geq \gamma \log_2 n $ , \\
 if $\tm>0,\
S_{\tM'}-\max_{0 \leq t \leq
 \tm}\left(S_t\right) \geq \gamma \log_2 n $ .
\end{enumerate}
%ou $\ki \equiv \ki(N)$ est tel que: $\lim_{N \rightarrow \infty}
%\ki(N)=0$ and $\lim_{N \rightarrow \infty} \ki(N)\ln N=+\infty$
\end{Def}

\noindent \textbf{The basic valley $\{{M_n}',\tmo,{M_n}\}$}

 We recall the notion of \textit{ basic valley } introduced
by Sinai and denoted here $\{{M_n}',\tmo,{M_n}\}$. The definition
we give is inspired by the work of \cite{Kesten2}. First let
$\{\tM',\tmo,\tM''\}$ be the smallest \textit{valley that contains
$0$ and of depth larger than} $\Gamma_n$. Here smallest means that
if we construct, with the operation of refinement, other valleys
in $\{\tM',\tmo,\tM''\}$ such valleys will not satisfy one of  the
properties of Definition \ref{thdefval1b}. ${M_n}'$ and ${M_n}$
are defined from $\tmo$ in the following way: if $\tmo>0$
\begin{eqnarray}
& & {M_n}'=\sup \left\{l\in \Z_-,\ l<\tmo,\ S_l-S_{\tmo}\geq
\Gamma_n,\
S_{l}-\max_{0 \leq k \leq \tmo}S_k \geq \gamma \log_2 n \right\} ,\\
& & {M_n}=\inf \left\{l\in \Z_+,\ l>\tmo,\ S_l-S_{\tmo}\geq
\Gamma_n\right\} . \label{4.8}
\end{eqnarray}
if $\tmo<0$
\begin{eqnarray}
& & {M_n}'=\sup \left\{l\in \Z_-,\ l<\tmo,\ S_l-S_{\tmo}\geq
\Gamma_n\right\} , \\
& & {M_n}=\inf \left\{l\in \Z_+,\ l>\tmo,\ S_l-S_{\tmo}\geq
\Gamma_n,\ S_{l}-\max_{ \tmo \leq k \leq 0}S_k \geq \gamma \log_2
n \right\} . \label{4.10}
\end{eqnarray}
if $\tmo=0$
\begin{eqnarray}
& & {M_n}'=\sup \left\{l\in \Z_-,\ l<0,\ S_l-S_{\tmo}\geq
\Gamma_n \right\} , \\
& &  {M_n}=\inf \left\{l\in \Z_+,\ l>0,\ S_l-S_{\tmo}\geq \Gamma_n
\right\} . \label{4.12}
\end{eqnarray}
\noindent  $\{{M_n}',\tmo,{M_n}\}$ exists
 with a $Q$ probability as close to one as we need. In fact it is not
 difficult to prove the following lemma

\begin{figure}[h]
\begin{center}
%\input{thfig4.pstex_t} \caption{} \label{thfignew1}
%\end{center}
%\end{figure}
%\begin{figure}[h]
%\begin{center}
\input{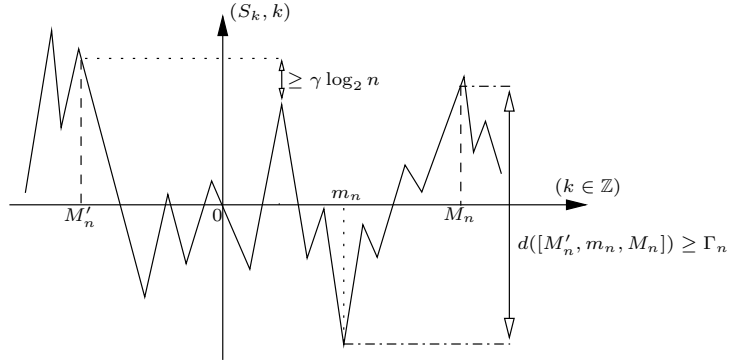}
 \caption{Basic valley, case $m_{n}>0$ } \label{thfig4}
\end{center}
\end{figure}

\begin{Lem} \label{moexiste} There exists $c>0$ such that if  \ref{hyp1}, \ref{hyp0} and \ref{hyp4} hold, for all $\gamma>0$ and $n$ we have
\begin{eqnarray}
Q\left[\{{M_n}',\tmo,{M_n}\} \neq \varnothing \right] =
1-\frac{c \gamma \log_{2} n}{\log n}.
\end{eqnarray}
\end{Lem}

\begin{Pre}
One can find the proof of this Lemma in Section 5.2 of \cite{Pierre2}.
\end{Pre}

\subsection{Main results}
We start with some definitions that will be used all along this work.
\noindent Let $x \in \Z$, define
\begin{eqnarray}
& &  T_{x}=\left\{\begin{array}{l} \inf\{k\in\N^*,\ X_k=x \}
\\ + \infty \textrm{, if such }  k \textrm{ does not exist.}
 \end{array} \right. \label{3.7sec}
\end{eqnarray}
Let $n>1$, $k \in \Z$, and $c_{0}>0$, define: 
\begin{eqnarray}
&& S_{k,m_{n}}^{n}=1-\frac{1}{\log n}(S_{k}-S_{m_{n}}), \\
&& \hat{S}_{k}^{n}=\frac{\log(\lo(k,n))}{\log n}, \\
% && u_{k,n}=  \left\{\begin{array}{ll}  \frac{c_{0}}{\log n} \textrm{ if } |k-m_{n}| \geq \log_{2}n, 
%\\ c_{0}'  \textrm{ if }   |k-m_{n}| < \log_{2}n.
&& u_{n}=\frac{c_{0}\log_{3} n}{\log n}.
% \end{array}
\end{eqnarray}
$S_{k,m_{n}}^{n}$ is the function of the potential we want to estimate,  $ \hat{S}_{k}^{n}$ is the estimator and $u_{n}$ is an error function.

\noindent Now let us define the following random sub-set of $\Z$:
\begin{eqnarray}
%\Ln^{\gamma}=\left\{k\in \Z,\ \sum_{j=T_{k^*}}^{n}\un_{X_{j}=k} \geq (\log %n)^{\gamma},\ \lo(k,T_{k^{*}}) <  \sum_{j=T_{k^*}}^{n}  \un_{X_{j}=k} \right\}
\Ln^{\gamma}=\left\{l\in \Z,\ \sum_{j=T_{k^*}}^{n}  \un_{X_{j}=l} \geq (\log n)^{\gamma}   \right\},
\end{eqnarray}
recall that $\gamma>0$. This set $\Ln^{\gamma}$ is fundamental for our result, we notice that it  depends only on the trajectory of the walk and more especially of its local time: $\Ln^{\gamma}$ is the set of points for which we are able to give an estimator of the the random potential. We will see that this set is large and contains a great amount  of the points visited by the walk (see Proposition \ref{prop3b}). 
We recall that $T_{k^*}$ is the first time the walk hit a favorite site.  In words, $l \in \Ln^{\gamma}$, if and only if:  The local time of the random walker in $l$ after the instant $T_k^*$ is large enough (larger than $(\log n)^{\gamma}$). 
Our main result is the following:

\begin{The} \label{th1} Assume \ref{hyp1bb},
\ref{hyp0} and \ref{hyp4} hold, there exists three constants $c_{0}$, $c_{1}$, $c_{2}$ and $c_2'$ such that for all $\gamma
> 6$, there exists $n_0$ such that for all $n>n_0$ there
exists $G_n\subset \Omega_1$ with $Q\left[G_n\right] \geq 1-\phi_{1}(n) $ and
\begin{eqnarray}
& & \inf_{\alpha \in G_{n}} \pa\left[ \bigcap_{k\in \Ln^{\gamma} } \left\{ \left | \hat{S}_{k}^n-S_{k,m_{n}}^n  \right |<u_{n} \right\}
 \right]  \geq 1-\phi_{2}(n).  \label{1.28}
\end{eqnarray}
where 
\begin{eqnarray}
& & \phi_{1}(n)= \frac{c_{1} \gamma \log_{2} n}{\log n}, \label{phi1} \\ 
& & \phi_{2}(n)=\frac{c_{2} }{(\log n)^{\gamma/2}}+\frac{c_{2} '}{(\log n)^{\gamma-6}} \label{phi2}. 
\end{eqnarray}
\end{The}

\noindent \\ The fact that our result depends on $m_{n}$ seems to be restrictive, we would like to know where is the bottom of the valley only by considering the local time of the walk, we prove the following:

\begin{Pro} \label{prop1}
Assume \ref{hyp1bb},
\ref{hyp0} and \ref{hyp4} hold, there exists a constant $c_{3}>0$ such that for all $\gamma
> 6$, there exists $n_0$ such that for all $n>n_0$ there exists $G_n\subset \Omega_1$ with $Q\left[G_n\right] \geq 1-\phi_{1}(n) $ and
\begin{eqnarray}
& & \inf_{\alpha \in G_{n}} \pao\left[ \max_{x\in \F_{n}}\left | m_{n}-x  \right | \leq (\log_{2} n)^2 \right]  \geq 1-\phi_{3}(n), \label{1P1.10} \\
& & \inf_{\alpha \in G_{n}}   \pao\left[ \left| T_{m_{n}}-T_{k^{*}} \right| \leq ( \log n)^3   \right] \geq  1- \phi_{3}(n),   \label{2P1.10} 
\end{eqnarray}
where $\phi_{3}(n)=c_{3}/(\log n)^{\gamma-6}$.
\end{Pro}

Notice that the distance between $m_n$ (coordinate of the point visited by the walk where the minimum of the potential is reached) and a favorite site is negligible comparing to a typical fluctuation of the walk (of order $(\log n)^2$). Thanks to Proposition \ref{prop1} we can replace \ref{1.28} in  Theorem \ref{th1} by
\begin{eqnarray}
& & \inf_{\alpha \in G_{n}} \pa\left[ \bigcap_{k\in \Ln^{\gamma} } \left\{ \left | \hat{S}_{k}^n-S_{k,k^*}^n  \right |<u_{n} \right\}
 \right]  \geq 1-\phi_{2}(n).  
 \end{eqnarray}

\noindent Now let us give a result giving the main properties of $\Ln^\gamma $.

\begin{Pro} \label{prop3b}
Assume \ref{hyp1bb} and
\ref{hyp0} and \ref{hyp4} hold, there exists a constant $c_{3}>0$ such that for all $\gamma
> 6$, there exists $n_0$ such that for all $n>n_0$ there exists $G_n\subset \Omega_1$ with $Q\left[G_n\right] \geq 1-\phi_{1}(n) $,
\begin{eqnarray}
& & \inf_{\alpha \in G_{n}} \pao\left[\lo(\Ln^{\gamma},n) =n(1-o(1)) \right]  \geq 1-\phi_{2}(n), \label{eq33} \\
& & \inf_{\alpha \in G_{n}}   \pao\left[ \left| \Ln^{\gamma} \right| Ê\thickapprox ( \log n)^2   \right] \geq 1-\phi_{2}(n),  \label{eq34} \\
%& & \inf_{\alpha \in G_{n}}   \pao\left[ \F_n \subset \Ln^{\gamma}  \right] \geq 1-\phi_{2}(n).    
%& & \inf_{\alpha \in G_{n}}   \pao\left[  \Ln^{\gamma} \textrm{ is connex exept may be at the border}   \right] \geq 1-\phi_{2}(n).    \label{eq35bc}
\end{eqnarray}
where $\lim_{n \rightarrow + \infty}o(1)=0$. 
\end{Pro}

\begin{Rem}  By definition we have $\F_n \subseteq \Ln^{\gamma}$.
% moreover if we define $\tilde{L}_n:=\{k \in \Z, \lo(k,n) >0\}$. 
%in probability $|\Ln^{\gamma}|/\tilde{L}_n=1-o(1)$ for $n$ large enough.
\end{Rem}

%\begin{Rem} About the variance of the estimator.
%This last fact means that we can only have a sharp result after the walk reach the sub-valley where it will spend the most of it's time. The link between $\Ln^{\gamma}$ and the random environment
%. Notice also that it is easy to prove that $\Ln^{\gamma}$ is large: of the order of the typical fluctuation of the walk $(\log n)^2$.
%There is another important remark we need to state now,  even the result is sharp, we were not able to have a good approximation of the potential in a small neighborhood of $m_{n}$,  especially for points at distance smaller than $(\log_{2} n)^2$ from $\tmo$. 
%\begin{center}
%--------------------------------------
%\end{center}
%\end{Rem}

\noindent
Theorem \ref{th1} is known to be the quenched result that means for a fixed environment $\alpha$, a simple consequence (see Remark \ref{decP}) is the following annealed result:

\begin{Cor} \label{th2}
Assume \ref{hyp1bb},
\ref{hyp0} and \ref{hyp4} hold, there exists three constants $c_{0}$, $c_{1}$ and $c_{2}$ such that for all $\gamma
> 6$, there exists $n_0$ such that for all $n>n_0$
\begin{eqnarray}
& & \p\left[ \bigcap_{k\in \Ln^{\gamma} } \left\{ \left | \hat{S}_{k}^{k}-S_{k,k^*}^n  \right |<u_{n} \right\}
 \right]  \geq 1-\phi(n),
\end{eqnarray}
where $\phi(n)=\phi_{1}(n)+\phi_{2}(n)$.
\end{Cor}
%\noindent
%It is important to notice that $\Ln^{\gamma}$ is a large set, indeed ...
We would like to notice, that for our purpose the result above is not very interesting, because the aim is to reconstruct one environment whereas the result above give the mean of the probability for the walk over all the possible environments.

\noindent \\ This paper is organized as follows. In Section 2 we
give the proof of Theorems \ref{th1} (we easily get the corrolary from Remark \ref{decP}), we have split this proof into two parts, the first one deals with the random environment and the other one with the random walk itself. In section 3 we sketch the proofs of Propositions \ref{prop1} and \ref{prop3b}. In Section 4, as an application of our result, we present an algorithm and some numerical simulations. For completeness, we recall in the appendix, some basic facts on birth and death processes.
% \tableofcontents
\section{Proof of Theorem \ref{th1}}

The proof of a result with a random environment involves both arguments and properties for the random environment and arguments for the random walk itself. I will start to give the properties I need for the random environment. Then we will use it to get the result for the walk. 

\subsection{Properties needed for the random environment}

\subsubsection{ Construction of ($G_{n}, n\in \N$) \label{sec11}}

Let $k$ and $l$ be in $\Z$, define 
\begin{eqnarray}
\sa{k}{l}=\Ea_{k}\left[\lo(l,T_{k})\right]
\end{eqnarray}
in the same way, let $A\subset \Z$, define 
\begin{eqnarray}
\sa{k}{A}=\sum_{l\in A}\Ea_{k}\left[\lo(l,T_{k})\right].
\end{eqnarray}

\begin{Def} \label{superb}  Let $d_{0}>0$, $d_{1}>0$, and $\omega \in \Omega_1$, we will say that $\alpha \equiv
\alpha(\omega)$ is a \textit{good environment} if there exists
$n_0$ such that for all $n\geq n_0$ the sequence $(\alpha_i,\ i
\in \Z)=(\alpha_i(\omega),\ i \in \Z)$ satisfies
 the properties \ref{3eq325} to \ref{8eq36}
\begin{eqnarray}
& \bullet &  \{{M_n}',\tmo, {M_n}\} \neq \varnothing,   \label{3eq325}\\
 & \bullet  & {M_n}'\geq -d_{0}( \sigma^{-1} \log_{2}n \log n)^2
 ,\ {M_n}\leq d_{0}( \sigma^{-1} \log_{2}n \log n)^2  ,
\label{interminibb} \\
& \bullet & \sa{m_{n}}{W_{n}} \leq  d_{1}  (\log_{2} n)^2, \label{8eq36}
\end{eqnarray}
where $W_{n}=\{M_{n}',M_{n}'+1,\cdots,m_{n},\cdots,M_{n}\}$.
%\begin{eqnarray}
%& \bullet & \textrm{basic properties to get the results of \cite{Pierre1}} \{}
%\end{eqnarray}
\end{Def}

\begin{Rem}
We will see in Section 2 that we use some results of \cite{Pierre2}. Considering this, we need extra properties on the random environment in addition to the three mentioned above, but as we don't need them for our computations we do not make them appear.
\end{Rem}

%\textbf{ Rapport du referee (minor comments 19) : appuyer sur le
%fait que la notion de bon environnement dépend de $n$}

\noindent  Define the \textit{set of good environments}
\begin{eqnarray}
G_n \equiv G_{n}(d_{0},d_{1})=\left\{\omega \in \Omega_1,\ \alpha(\omega) \textrm{ is a }
\textit{ good environment} \right\}.
\end{eqnarray}
$G_n$ depends on $d_0$, $d_{1}$ and $n$, however  we only make explicit the $n$ dependence.

\begin{Pro} \label{profondab} There exists two constants $d_{0}>0$ and $d_{1}>0$ such that if   \ref{hyp1bb}, \ref{hyp0} and \ref{hyp4} hold,  there exists $n_0$ 
such that for $n>n_0$
\begin{eqnarray}
Q\left[ G_n\right]  \geq 1- \phi_{1}(n),
\end{eqnarray}
where $\phi_{1}(n)$ is given by \ref{phi1}.
\end{Pro}

\begin{Pre}
We can find the proof for the first three properties \ref{3eq325}-\ref{8eq36} in \cite{Pierre2}, see Definition 4.1 and Proposition 4.2. 
\end{Pre}

\noindent \\ To end the section we would like to make the following elementary remark on the decomposition of $\p$:

\begin{Rem} \label{decP}
Let $\Ct_{n} \in \sigma \left(X_{i}, i \leq n \right) $ and $G_{n} \subset \Omega_{1}$, we have :
\begin{eqnarray}
\p\left[\Ct_{n}\right] &\equiv & \int_{\Omega_{1}}Q(d\omega)\int_{\Ct_{n}} d\p^{\alpha(\omega)} \\
 & \geq &  \int_{G_{n}}Q(d\omega)\int_{\Ct_{n}} d\p^{\alpha(\omega)}.
\end{eqnarray}
So assume that $Q[G_{n}] \equiv e_{1}(n) \geq 1-\phi_{1}(n)$ and assume that for all $\omega \in G_{n}$, $\int_{\Ct_{n}} d\p^{\alpha(\omega)} \equiv e_{2}(\omega,n) \geq 1-\phi_{2}(n) $ we get that
\begin{eqnarray}
\p\left[\Ct_{n}\right] \geq e_{1}(n) \times \min_{w\in G_{n}}(e_{2}(w,n)) \geq 1-\phi_{1}(n)-\phi_{2}(n).
\end{eqnarray}
\end{Rem}
%So choosing $\Ct_{n}= \left\{\max_{x}\lo\left(\teb,n\right) \geq \beta n \right\}$,  we have to extract from $\Omega_{1}$ a subset $G_{n}$ sufficiently small to get that $\min_{w\in G_{n}}(d_{1}(w,n))>0$ (Proposition \ref{Prop2}) but sufficiently large to have $d_{2}(n)>0$ (Proposition \ref{profondab}) .
%The largest part of the proof is to construct such a $G_{n}$ (Section \ref{sec11} and Appendix B).

\subsection{Arguments for the walk}

Let $(\rho_{1}(n),n\in \N)$ a strictly positive decreasing sequence such that $\lim_{n \rightarrow \infty} \rho_{1}(n)=0 $. First let us show that the Theorem \ref{th1} is a simple consequence of the following

\begin{Pro} \label{ProF1} 
Assume \ref{hyp1bb},
\ref{hyp0} and \ref{hyp4} hold, there exists $n_0$ such that for all $n>n_0$ there
exists $G_n\subset \Omega_1$ with $Q\left[G_n\right] \geq 1-\phi_{1}(n) $ and
\begin{eqnarray}
 & & \quad  \sup_{\alpha \in G_n}
\left\{\pao\left[\bigcup_{k\in
\Ln^{\gamma}}\left\{\left|\frac{\lo(k,n)}{n}-
\frac{\sa{m_{n}}{k}}{\sa{m_{n}}{W_{n}}} \right| \leq w_{k,n} \right\} \right] \right\} \geq 1- \phi_{2}(n) \label{extension}
\end{eqnarray}
where $ w_{k,n}=\rho_{1}(n)\frac{\sa{m_{n}}{k}}{\sa{m_{n}}{W_{n}}} $, $\phi_{1}(n)$ and $\phi_{2}(n)$ are given just after \ref{1.28}
\end{Pro}
\noindent
Taking the logarithm and for $n$ large enough, using Taylor series expansion, we remark that 
\begin{eqnarray}
\frac{\sa{m_{n}}{k}}{\sa{m_{n}}{W_{n}}}(1-\rho_{1}(n))  \leq \frac{\lo(k,n)}{n}
\leq \frac{\sa{m_{n}}{k}}{\sa{m_{n}}{W_{n}}} \left(1+ \rho_{1}(n)\right)
\end{eqnarray}
implies
\begin{eqnarray}
& & -2\rho_{1}(n)-\log(\sa{m_{n}}{W_{n}})   \leq \log \lo(k,n)-\log n-  \log (\sa{m_{n}}{k})  \leq -\log(\sa{m_{n}}{W_{n}}) + \rho_{1}(n), \nonumber
\end{eqnarray}
rearranging the terms and using \ref{A.1} (see the Appendix) we get
\begin{eqnarray}
& & \frac{1}{\log n}(\ra_{n}(k)-2\rho_{1}(n) )  \leq  \hat{S}_{k}^n-S_{k,m_{n}}^n  \leq \frac{1}{\log n}(\ra_{n}(k)- \rho_{1}(n))
\end{eqnarray}
where $\ra_{n}(k)=\log\left( \frac{\alpha_{\tmo}}{\beta_k}  a_{k,m_{n}}\right )-\log(\sa{m_{n}}{W_{n}}) $ and $ a_{k,m_{n}}$ is given by \ref{A.2}. Now using \ref{A.4} and Property \ref{8eq36} we get the Theorem.
\noindent
The proof of Proposition \ref{ProF1} is based on the following results (Lemma \ref{LemA1}) of \cite{Pierre1}, 

\subsubsection{Known facts}

Let $(\rho(n),n \in \N)$ be a positive decreasing sequence such that $\lim_{n\rightarrow \infty} \rho(n)=0$, we define
\begin{eqnarray}
&  & \A_{1}=\left\{ \left|\frac{\lo(\tmo,n)}{n}-
\frac{1}{\sa{m_{n}}{W_n}} \right|> \frac{\rho(n)}{\sa{m_{n}}{W_n}} \right\}, \\
& & \A_{2}=\left\{ T_{m_{n}}\leq n/(\log n)^4, \lo(W_{n},n)=1 \right\}.
\end{eqnarray}

\begin{Lem} \label{LemA1} 
Assume \ref{hyp1bb},
\ref{hyp0} and \ref{hyp4} hold, there exists a constant $b_{1}>0$ such that for all $\gamma
>6$, there exists $n_0$ such that for all $n>n_0$ there
exists $G_n\subset \Omega_1$ with $Q\left[G_n\right] \geq 1-\phi_{1}(n) $ and
\begin{eqnarray}
 & \quad &  \sup_{\alpha \in G_n'}
\left\{\pao\left[\A_{1}\right] \right\} \leq r_{1}(n), \label{10eq43}
\end{eqnarray}
where $r_{1}(n)= b_{1}/(\log n)^{\gamma-6}$.
\end{Lem}

\begin{Pre}
We do not give the details of the computations because the reader can find it in the referenced paper (Theorem 3.8 of \cite{Pierre2}), just notice that comparing to the Theorem 3.8 we have a better  rate of convergence for the probability obtained just by using a weaker result for the concentration of the walk.
\end{Pre}

\noindent \\ We will also need the following elementary fact :
\begin{Lem} \label{LemA2} 
Assume \ref{hyp1bb},
\ref{hyp0} and \ref{hyp4} hold, there exists a constant $b_{2}>0$ such that for all $\gamma
>2$, there exists $n_0$ such that for all $n>n_0$ there
exists $G_n\subset \Omega_1$ with $Q\left[G_n\right] \geq 1-\phi_{1}(n) $ and \begin{eqnarray}
 & \quad &  \sup_{\alpha \in G_n'}
\left\{\pao\left[\A_{2}  \right] \right\} \leq r_{2}(n), \label{10eq43}
\end{eqnarray}
where $r_{2}(n)=  b_{2}/(\log n)^{\gamma-2}$. 
\end{Lem}

\begin{Pre}
Once again this can be find in  \cite{Pierre2}: Proposition 4.7 and Lemma 4.8.
\end{Pre}

\noindent \\
Using these results we can give the proof of Proposition  \ref{ProF1} into two steps :

\subsubsection {Step 1}
\begin{figure}[h]
\begin{center}
%\input{thfig4.pstex_t} \caption{} \label{thfignew1}
%\end{center}
%\end{figure}
%\begin{figure}[h]
%\begin{center}
\input{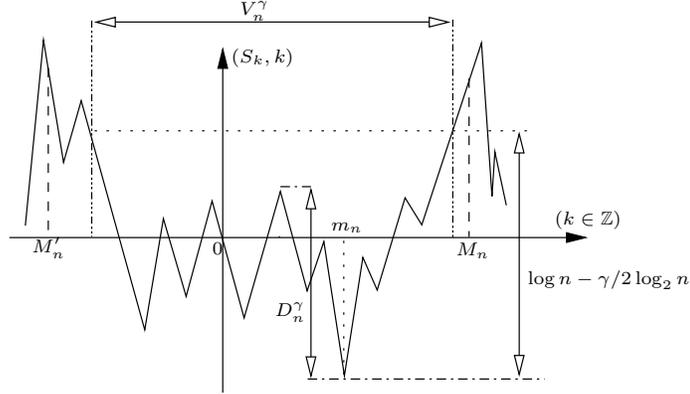}
 \caption{ $V_{n}^{\gamma}$ with $m_{n}>0$, case  1: $ D_{n}^{\gamma} \equiv  \max_{ k \in \bar{v}_1} \max_{k \leq j \leq m_{n}}(S_{j}-S_{m_{n}}) \leq \log n - \gamma \log_{2} n$} \label{fig6}
\end{center}
\end{figure}

Let us define the following subsets :
\begin{eqnarray}
 \bar{v}_1^n & \equiv & \{M_n' \leq k \leq m_n-1,\ ( \max_{k \leq j \leq m_n}S_{j}-S_{m_{n}})< \log n -\frac{\gamma}{2} \log_{2} n \},  \\
 \bar{v}_2^n & \equiv & \{m_n+1 \leq k \leq M_n,\ ( \max_{m_n \leq j \leq k}S_{j}-S_{m_{n}})< \log n - \frac{\gamma}{2} \log_{2} n \}, 
\end{eqnarray}
and 
\begin{eqnarray}
& & V_{n}^{\gamma}=  \bar{v}_1^n  \cap  \bar{v}_2^n . \label{DefVn}
\end{eqnarray}
In words $ V_{n}^{\gamma}$ is a subset of points included in $W_n$,  such that for all $k \in V_{n}^{\gamma}$ the largest difference of potential between $m_n$ and $k$ is smaller than $\log n -\gamma/2 \log_{2} n$. For the walk we will see (Lemma below) that if $k \in V_{n}^{\gamma}$ then the walk will hit $k$ after it has reached $m_n$ and it will hit this point $k$ a number of time large enough (see figure \ref{fig6}).

%\begin{figure}[h]
%\begin{center}
%\input{fig6b.pstex_t}
% \caption{ $V_{n}^{\gamma}$ with $m_{n}>0$, case  2: $D_{n}^{\gamma} \equiv \max_{ k \in \bar{v}_1} \max_{k \leq j \leq m_{n}}(S_{j}-S_{m_{n}}) > \log n - \gamma \log_{2} n$} \label{fig6b}
%\end{center}
%\end{figure}

First let us prove the following Lemma :
\begin{Lem} \label{LemA3} 
Assume \ref{hyp1bb},
\ref{hyp0} and \ref{hyp4} hold, there exists a constant $b_{3}>0$ such that for all $\gamma
>6$, there exists $n_0$ such that for all $n>n_0$ there
exists $G_n\subset \Omega_1$ with $Q\left[G_n\right] \geq 1-\phi_{1}(n) $\begin{eqnarray}
 & \quad &  \sup_{\alpha \in G_n'}
\left\{\pao\left[ \Ln^{\gamma} \subseteq V_{n}^{\gamma}  \right] \right\} \geq  1-r_{3}(n), \label{10eq43}
\end{eqnarray}
where $r_{3}(n)= b_{3}/ (\log n)^{\gamma / 2}  $. 
\end{Lem}

\noindent Notice that $\Ln^{\gamma}$ is a $\p$ random variable (with two levels of randomness) whereas $V_{n}^{\gamma}$ is only a $Q$ random variable (with one level of randomness), this Lemma makes the link between a trajectory of the walk and the random environment.  %$ \Ln^{\gamma} \subset V_{n}^{\gamma}$ informs us what are the property of the environment to get $\Ln^{\gamma}$.

\begin{Pre}
To prove this Lemma we use Proposition \ref{prop1}. First notice that \begin{eqnarray}
\pao\left[ \Ln^{\gamma} \subseteq V_{n}^{\gamma}  \right]  = 1-\pao\left[ \bigcup_{k\in (v_{n}^{1}\cup v_{n}^{2}) } \left\{ k\in \Ln^{\gamma}  \right\} \right]  
\end{eqnarray}
where 
 \begin{eqnarray}
 v_1^n & \equiv & \{M_n' \leq k \leq m_n-1,\ ( \max_{k \leq j \leq m_n}S_{j}-S_{m_{n}}) \geq  \log n -\frac{\gamma}{2} \log_{2} n \},  \\
 v_2^n & \equiv & \{m_n+1 \leq k \leq M_n,\ ( \max_{m_n \leq j \leq k}S_{j}-S_{m_{n}}) \geq \log n - \frac{\gamma}{2} \log_{2} n \}. 
   \end{eqnarray} 
\noindent 
Let $k \in v_{1}^{n}$ let us give an upper bound for  
\begin{eqnarray}
\qquad \pao\left[ k \in \Ln^{\gamma},\   |T_{k^*}-T_{m_{n}}| \leq (\log n)^3 \right] & \leq & \pao\left[\sum_{j=T_{k^*}}^{n}\un_{X_{j}=k} \geq (\log n)^{\gamma},\ |T_{k^*}-T_{m_{n}}| \leq (\log n)^3 \right] \nonumber \\
& \leq & \pao\left[\sum_{j=T_{m_{n}}}^{n}\un_{X_{j}=k} \geq (\log n)^{\gamma}- (\log n)^3 \right]  \nonumber \\
& \leq & \pam\left[\sum_{j=1}^{T_{m_{n},n}}\un_{X_{j}=k} \geq (\log n)^{\gamma}-(\log n)^3\ \right],  
\end{eqnarray}
for the third inequality we have used the strong Markov property, where
\begin{eqnarray*}
&& T_{\tmo,j } \equiv \left\{\begin{array}{l}  \inf\{k>T_{\tmo,j-1},\ X_k=\tmo \},\  j \geq 2  \\
 + \infty \textrm{, if such }  k \textrm{ does not exist.}
 \end{array} \right. \\
& & T_{\tmo,1}  \equiv T_{\tmo} \ (\textrm{see \ref{3.7sec}}).
\end{eqnarray*}
Now using the Markov inequality and Lemma \ref{Lou} we get
\begin{eqnarray}
 \pao\left[ k \in \Ln^{\gamma},\   |T_{k^*}-T_{m_{n}}| \leq (\log n)^3\right]  & \leq & \frac{n \Eam\left[\lo(k,T_{m_{n}})\right]}{(\log n)^{\gamma}-(\log n)^3} \\
& \leq & \frac{n}{\eta_{0}\exp(S_{k}-S_{m_{n}})((\log n)^{\gamma}-(\log n)^3)} \\
& \leq & \frac{1}{\eta_{0} (\log n)^{\gamma/2}(1-(\log n)^3/(\log n)^{\gamma})},
\end{eqnarray}
notice that in the last inequality we have used the fact that $k\in v_1^{n}$.  A similar computation give the same inequality when $kÊ\in v_2^n$ . Collecting what we did above, and using the Property \ref{interminibb} together with \ref{2P1.10} yields the Lemma.
 \end{Pre}

\subsubsection{Step 2}

\noindent This second step is devoted to the proof of the following Lemma.
%Now we prove the proposition, let us prove the following (quenched) Lemma, denoting $\A_{1}=
%\left|\frac{\lo(\tmo,n)}{n}-\frac{1}{\sa{m_{n}}{W_n}} \right| \leq \rho(n) $ and $\A_{2}=T_{m_{n}}\leq %0n^{\delta_{1}}, \lo(W_{n},n)=1$, we have :

\begin{Lem} \label{LemF2} 
For all $\alpha$ and $n$ we have 
\begin{eqnarray}
 & & \pao\left[
\left|\frac{\lo(k,n)}{n}-
\frac{\sa{m_{n}}{k}}{\sa{m_{n}}{W_{n}}} \right|> w_{k,n},\ \A_{1},\  \A_{2} \right]  \leq 2 \exp(-n/2\psi_{2}^{\alpha}(n)) \label{extension}
\end{eqnarray}
recall that  $ w_{k,n}=\rho_{1}(n)\frac{\sa{m_{n}}{k}}{\sa{m_{n}}{W_{n}}}  $ and $\psi_{2}^{\alpha}(n)= 2\frac{(\rho_{1}(n)-\rho(n))^{2}}{1+\rho(n)}  \frac{(\alpha_{m_{n}} \wedge \beta_{m_{n}})}{|k-m_{n}|}\frac{\exp(-(S_{M_{k}}-S_{m_{n}}))}{\sa{m_{n}}{W_{n}}}$. $M_{k}$ is such that $S_{M_{k}}=\max_{ m_{n}+1 \leq j \leq k } S_{j}$ if $k>m_{n}$ and conversly if $k<m_{n}$ $S_{M_{k}}=\max_{  k \leq j \leq m_{n}-1 } S_{j}$.
\end{Lem}

\begin{Pre}
We essentially use an inequality of concentration (see \cite{Ledoux}), for simplicity we only give the proof for $k>m_{n}$, the other case ($k \leq m_{n}$) is very similar. Using  the Markov property and the fact that $\lo(k,T_{m_{n}})=0$, we get 
\begin{eqnarray}
& & \pao\left[\left|\frac{\lo(k,n)}{n}-
\frac{\sa{m_{n}}{k}}{\sa{m_{n}}{W_{n}}} \right|> w_{k,n},\A_{1}, \A_{2} \right] \leq \pam\left[\left|\frac{\lo(k,n)}{n}-
\frac{\sa{m_{n}}{k}}{\sa{m_{n}}{W_{n}}} \right|> w_{k,n},\A_{1} \right]. 
\end{eqnarray}
We have 
\begin{eqnarray}
& & \pam\left[\frac{\lo(k,n)}{n}-
\frac{\sa{m_{n}}{k}}{\sa{m_{n}}{W_{n}}} > w_{k,n},\A_{1} \right] \nonumber \\
& \leq &  \pam\left[\frac{\lo(k,n)}{n}-
\frac{\sa{m_{n}}{k}}{\sa{m_{n}}{W_{n}}} > w_{k,n},  \frac{\lo(\tmo,n)}{n}-\frac{1}{\sa{m_{n}}{W_n}}  \leq \frac{\rho(n)}{\sa{m_{n}}{W_{n}}} \right] \\
&\leq & \pam\left[\frac{\lo(k,T_{m_{n},n_{1}})}{n}-
\frac{\sa{m_{n}}{k}}{\sa{m_{n}}{W_{n}}} > w_{k,n}\right]  \\
& \equiv &  \pam\left[\frac{\lo(k,T_{m_{n},n_{1}})}{n}-
\frac{\sa{m_{n}}{k}}{\sa{m_{n}}{W_{n}}}(1+\rho(n)) > w_{k,n}'\right]  
\end{eqnarray}
where $n_{1}=\frac{n}{\sa{m_{n}}{W_{n}}}\left(1+\rho(n) \right)$, notice that $n_{1}$ is not necessarily an integer but for simplicity we disregard that, and $ w_{k,n}'=\frac{\sa{m_{n}}{k}}{\sa{m_{n}}{W_{n}}}(\rho_{1}(n)-\rho(n)) $. The strong Markov property implies that $\lo(k,T_{m_{n},n_{1}})$ is a sum of $n_{1}$ i.i.d. random variables, the inequality of concentration gives 
 \begin{eqnarray}
& &   \pam\left[\frac{\lo(k,T_{m_{n},n_{1}})}{n}-
\frac{\sa{m_{n}}{k}}{\sa{m_{n}}{W_{n}}} > w_{k,n}',\A_{1} \right]  
\leq \exp\left[-\frac{n}{2} \frac{\sa{m_{n}}{W_{n}}}{\V_{m_{n}}(\lo(k,T_{m_{n}}))} \frac{ (w_{k,n}')^2}{1+\rho(n)} \right ].
\end{eqnarray}
With the same method we also get
 \begin{eqnarray}
& &  \pam\left[\frac{\lo(k,T_{m_{n},n_{1}})}{n}-
\frac{\sa{m_{n}}{k}}{\sa{m_{n}}{W_{n}}} <- w_{k,n}',\A_{1} \right]  
\leq \exp\left[-\frac{n}{2} \frac{\sa{m_{n}}{W_{n}}}{\V_{m_{n}}(\lo(k,T_{m_{n}}))} \frac{(w_{k,n}')^2}{1+\rho(n)} \right ].
\end{eqnarray}
Using \ref{A.3} we get Lemma \ref{LemF2}.
\end{Pre}

\subsubsection{End of the proof of the Theorem}

Using Lemmata \ref{LemA1}, \ref{LemA2} and \ref{LemA3} we have:

\begin{eqnarray}
  & &
\pao\left[\bigcup_{k\in
\Ln^{a}}\left\{\left|\frac{\lo(k,n)}{n}-
\frac{\sa{m_{n}}{k}}{\sa{m_{n}}{W_{n}}} \right|> w_{k,n} \right\} \right]  \nonumber \\
& \leq &\left|V_{n}^{\gamma} \right| \sup_{k\in V_{n}^{\gamma} } \pao\left[ \left|\frac{\lo(k,n)}{n}-
\frac{\sa{m_{n}}{k}}{\sa{m_{n}}{W_{n}}} \right|> w_{k,n},\A_{1}, \A_{2}  \right] +3\max_{1 \leq i\leq 3}\{r_{i}(n)\}, \nonumber
\end{eqnarray}
then using Lemma \ref{LemF2} we get 
\begin{eqnarray}
  \sup_{k\in V_{n}^{\gamma} } \pao\left[ \left|\frac{\lo(k,n)}{n}-
\frac{\sa{m_{n}}{k}}{\sa{m_{n}}{W_{n}}} \right|> w_{k,n},\A_{1}, \A_{2}  \right] & \leq &
 2 \sup_{k\in V_{n}^{\gamma} } \exp(-n/2 \psi_{2}^{\alpha}(k,n)) \nonumber \\
 & \leq & 2 \exp(- (\log n)^{\gamma/2-2}/(\rho_{1}(n) \log_{2} n)  ),
\end{eqnarray}
where the last inequality comes from the definition of $V_{n}^{\gamma} $ (see \ref{DefVn}) and the Properties \ref{interminibb} and \ref{8eq36}. To end we use again the Property \ref{interminibb} together with the definition of $V_{n}^{\gamma} $.

\section{Proof of Proposition \ref{prop1} and \ref{prop3b}}

\textit{Sketch of the proof of Proposition \ref{prop1}}
\ref{1P1.10} is an improvement of the proof of Corollary 3.17 of \cite{Pierre2}  in order to get a better rate of convergence for the probability. To get \ref{2P1.10}, we have used the same idea of the proof of Corollary 3.17 of \cite{Pierre2}, so once again we will not repeat the computations here. We recall just the intuitive idea: once the walk has reached $k^*$, we know from \ref{1P1.10} that $m_{n}$ is at most at a distance $(\log_{2}n)^{2}$, therefore the walk need at most the amount of time $\exp ( \sqrt((\log_{2}n)^{2}))=(\log n)$ to reach $m_n$. We take $(\log n)^3$ to get a better rate of convergence for the probability. \\
\textit{Sketch of the proof of Proposition \ref{prop3b}} 
The first two properties can be deduced from the following inequality, let $\epsilon>1$, for all $n$ large enough and all $\alpha \in G_n$:
\begin{eqnarray}
\pao\left[ V^{2(\gamma+\epsilon)}_n \subseteq \Ln^{\gamma}  \right] \geq \phi_3(n)+r_1(n)+cte(\log t)^2\exp\left(-(\log n)^{\gamma+\epsilon-2}(1-(\log n)^{1-\epsilon})\right). \label{eq10}
\end{eqnarray}
Indeed,thanks to \ref{eq10} we have 
\begin{eqnarray}
\pa\left[\lo(\Ln^{\gamma},n) \geq n(1-o(1)) \right] \geq \pa\left[\lo(V_n^{2(\gamma+ \epsilon)},n) \geq n(1-o(1)) \right] 
\end{eqnarray}
we get \ref{eq33} by using the same method \cite{Pierre2} uses to get Theorem 3.1.
To get \ref{eq34}, we only need to show that $|V_n^{\gamma+\epsilon}| \thickapprox (\log t)^2$, which is a basic fact for a simple random walk.
Now, to get \ref{eq10}, first we notice that, by using a similar method of the proof of Theorem \ref{th1} we can get
\begin{eqnarray}
\pao\left[ V^{2(\gamma+\epsilon)}_n \nsubseteq \Ln^{\gamma}  \right] \leq |V^{\gamma+\epsilon}_n| \max_{k \in V^{2(\gamma+\epsilon)}_n} \pam\left[\sum_{j=1}^{n_1}\eta_{j}^k< (\log n)^{\gamma} \right] + \phi_3(n)+r_1(n). 
\end{eqnarray}
where $(\eta_{j}^k,j)$ is a i.i.d. sequence with the law of $\lo(k,T_{m_n})$. Then using an inequality of concentration, we get \ref{eq10}.
%Let $\tilde{k}_1= \min\{j, j \in \Ln^{\gamma}\}$ and $\tilde{k}_2= \max\{j, j \in \Ln^{\gamma}\}$, $\bar{k}_1= \min\{j, j \in V_n^{\gamma}\}$ and $\bar{k}_2= \max\{j, j \in V_n^{\gamma}\}$
%\begin{eqnarray}
% \pao\left[ \bigcup_{\tilde{k}_1<k<\tilde{k}_2} \sum_{j=T_{k^*}}^{n}\un_{X_{j}=k} < (\log n)^{\gamma} \right] & \leq &  \pao\left[ \bigcup_{\bar{k}_1<k<\bar{k}_2} \sum_{j=T_{k^*}}^{n}\un_{X_{j}=k} < (\log n)^{\gamma} \right] +r_3(n) \\
 %& \equiv & \pao\left[\bigcup_{k \in V_n^{\gamma}} \sum_{j=T_{k^*}}^{n}\un_{X_{j}=k} < (\log n)^{\gamma} \right] +r_3(n),
%\end{eqnarray}
%the first inequality comes from \ref{10eq43}, last equality comes from the fact that $V_n^{\gamma}$ is connex by definition. So we get \ref{eq35bc}, using \ref{10eq43} once again. 

%Indeed it is proven that that all the favorite site belongs to a small interval centered in $\m_{n}$, this %give the first inequality of the Proposition. The second 

%\begin{eqnarray}
%  \pao\left[ \left| m_{n}-k^{*} \right|> f(n)  \right] \leq   \pao\left[ \left| m_{n}-k^{*} \right|> f(n), \A_{1},\A_%{2} \right]+ \pao\left[\A_{1} \right]+\pao\left[\A_{2} \right]
%\end{eqnarray}
%Let us give an upper bound for $\pao\left[ k^{*}-m_{n} > f(n),\ \A_{1},\ \A_{2}  \right] $, we have :
%\begin{eqnarray}
%\pao\left[ k^{*} -m_{n} > f(n),\ \A_{1},\ \A_{2}  \right]  \leq \pam\left[  \lo(m_{n},n_{1}) \geq \lo(m_{n}+f%(n),n)   \right] 
%\end{eqnarray}
%recall that $n_{1}=...$.

%\section {Proof of Proposition \ref{Prop2}}

\section{Algorithm and Numerical simulations}

\subsection{General and recall of the main definitions}
First notice that we have no criteria to determine wether or not we can apply this method to an unknown series of data. All we know is that it works for Sinai's walk, however we can apply the following algorithm to every process.
Let us recall the basic random variables that will be used for our simulations, let $x \in \Z$, $n \in \N$, 
\begin{eqnarray}
& &  T_{x}=\left\{\begin{array}{l} \inf\{k\in\N^*,\ X_k=x \}
\\ + \infty \textrm{, if such }  k \textrm{ does not exist.}
 \end{array}, \right. \\
& &  \lo\left(x,n\right) \equiv \sum_{i=1}^n \un_{\{X_i=x\}} , 
\end{eqnarray}
\begin{eqnarray}
& & \lo^*(n)=\max_{k\in \Z}\left(\lo(k,n)\right) ,\
\F_{n}=\left\{k \in \Z,\ \lo(k,n)=\lo^*(n) \right\} , \\
& & k^*=\inf\{|k|, k\in \F_{n} \}. \\
\end{eqnarray}
We recall also the set  $\Ln^{\gamma}$, the function of the potential we want to estimate and its estimator:
\begin{eqnarray}
& & \Ln^{\gamma}=\left\{k\in \Z,\ \sum_{j=T_{k^*}}^{n}  \un_{X_{j}=k} \geq (\log n)^{\gamma}   \right\}, \\
&& S_{k,m_{n}}^{n}=1-\frac{1}{\log n}(S_{k}-S_{m_{n}}), \\
&& \hat{S}_{k}^{n}=\frac{\log(\lo(k,n))}{\log n},
% && u_{k,n}=  \left\{\begin{array}{ll}  \frac{c_{0}}{\log n} \textrm{ if } |k-m_{n}| \geq \log_{2}n, 
%\\ c_{0}'  \textrm{ if }   |k-m_{n}| < \log_{2}n.
%&& u_{k,n}=\frac{c_{0}\log_{3} n}{\log n}.
% \end{array}
\end{eqnarray}
We also recall that thanks to Proposition \ref{prop1}, in probability we have $|m_n -k^*| \leq \textrm{cte} (\log_2 n)^2$. \\
\subsection{Main steps of the algorithm}
Step 1: We have to determine $\Ln^{\gamma}$ and to get it we have to compute $T_{k^*}$ and therefore the local time of the process. First we compute $\lo(k,n)$ for every $k$, notice that $\lo(k,n)$ is not equal to zero only if $k$ has been visited by the walk within the interval of time $[1,n]$. Then we can compute $\lo^*(n)$ and determine $k^*$ and $T_{k^*}$. Notice that $T_{k^*}$ is not a stopping time, therefore we need two passes to compute what we need. We are now able to determine $\Ln^{\gamma}$ computing $\sum_{j=T_{k^*}}^{n}  \un_{X_{j}=k}$. \\
Step 2: We can check that $\Ln^{\gamma}$ is connex, contains $k^*$ and that its size is of the order of a typical fluctuation of the walk. Now, keeping only the $k$ that belongs to $\Ln^{\gamma}$ we compute for those $k$: $ \hat{S}_{k}^{n}=\frac{\log(\lo(k,n))}{\log n}$ the estimator of the potential. We localize the bottom of the valley $m_n$ using $k^*$.
\subsection{Simulations}
\begin{figure}[h]
\begin{center}
\includegraphics*[width=11cm,height=5cm]{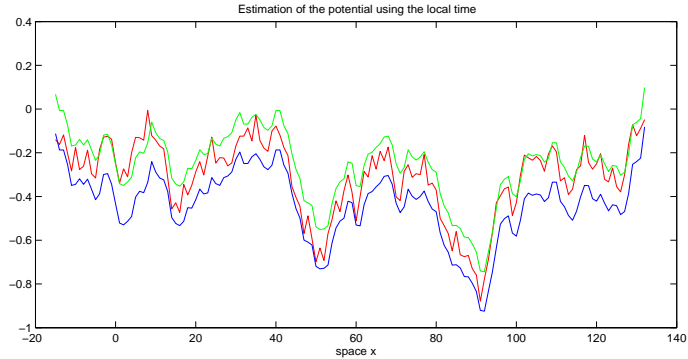}
\caption{in red $S_{x,m_{n}}^{n}$, in blue $ \hat{S}_{x}^{n}-u_n$, in green $ \hat{S}_{x}^{n}+u_{n}$} \label{fig7}
\end{center}
\end{figure}

\begin{figure}[h]
\begin{center}
\includegraphics*[width=11cm,height=5cm]{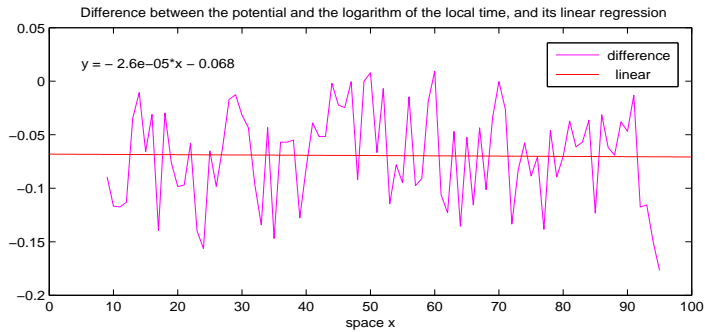}
\caption{in magenta $S_{x,m_{n}}^{n}-\hat{S}_{x}^{n}$, in red the linear regression} \label{fig8}
\end{center}
\end{figure}

\begin{figure}[h]
\begin{center}
\includegraphics*[width=11cm,height=5cm]{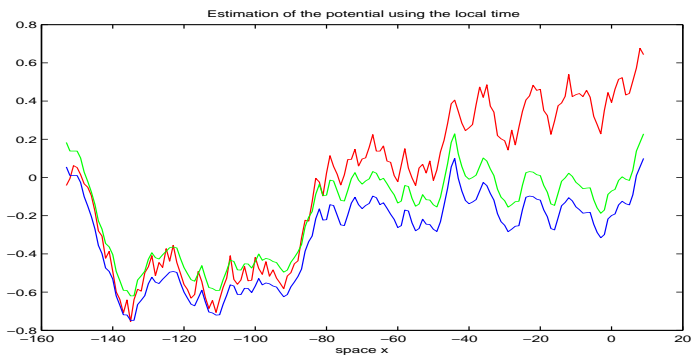}
\caption{in red $S_{x,m_{n}}^{n}$, in blue $ \hat{S}_{x}^{n}-u_n$, in green $ \hat{S}_{x}^{n}+u_{n}$} \label{fig9}
\end{center}
\end{figure}

\begin{figure}[hh]
\begin{center}
\includegraphics*[width=11cm,height=5cm]{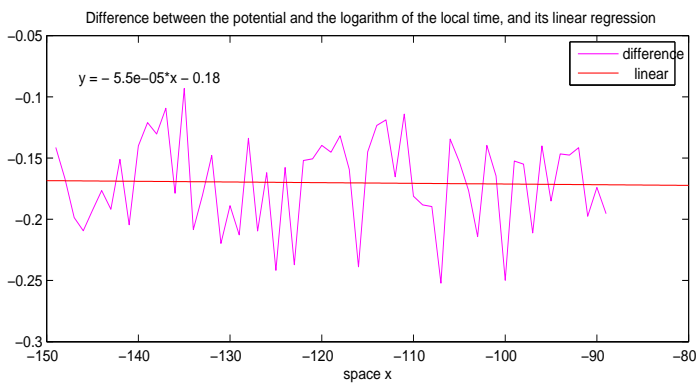}
\caption{in magenta $S_{x,m_{n}}^{n}-\hat{S}_{x}^{n}$, in red the linear regression} \label{fig10}
\end{center}
\end{figure}

For the first simulation (Figure \ref{fig7}) we show a case where $\Ln^{\gamma}$ is large i.e. $\Ln^{\gamma}$ contains most of the points visited by the walk. 
The trajectory of the random potential is in red the interval of confidence in blue and green. We took $n=500000$ and $\gamma = 7$, notice that the larger is $\gamma$, the  smaller is $\Ln^{\gamma}$ but better is the rate of convergence of the probability.
We get that $\Ln^{\gamma}=[10,94]$. In Figure \ref{fig8} we plot the difference $S_{x,m_{n}}^{n}-\hat{S}_{x}^{n}$ and its the linear regression.
We notice that the slope of the linear regression  is of order $10^{-5}$. We also notice that we have taken $n=500000$, so the error function $u_n  \thickapprox \frac{\log_3 n}{\log n} \thickapprox 0,7 $ this match with the $\max_{x}(S_{x,m_{n}}^{n}-\hat{S}_{x}^{n})\thickapprox 0.8  $ for this simulation.
\newpage
Now let us choose another example where $\Ln^{\gamma}$ is much more smaller. For the following simulation (Figure \ref{fig9}) we have only changed the sequence of random number.  
We get that $\Ln^{\gamma}=[-150,-85]$. We notice that for the coordinates larger than -85 and especially  after -40, our estimator is not good at all. In fact once the walk has reached the minimum of the valley (coordinate -111) it will never reach again one of the points of coordinate larger than  -40 before $n=500000$, so our estimator can not say anything about the difference $S_{x,m_{n}}^{n}-\hat{S}_{x}^{n}$. However  if we look in the past of the walk and especially at a the time $T_{k^*}$ which is the first time it has reached the coordinate $-111$, the favorite point for this time is localized around the point $-2$, so a good estimator between the coordinate -40 and 10 may be given by  $(\frac{\log(\lo(k,T^*))}{\log T^*},k)$.
The difference $S_{x,m_{n}}^{n}-\hat{S}_{x}^{n}$ and the  linear regression in the interval $\Ln^{\gamma}=[-150,-85]$ is presented Figure \ref{fig10}.
%\newpage
\newpage
\appendix

\section{Basic results for birth and death processes}

\noindent 
%\cite{Chung} and \cite{Revesz} on inhomogeneous discrete time
%birth and death processes.

%\noindent Let $x,a$ and $b$ in $\Z$, assume $a<x<b$,
%the two following lemmata can be found in \cite{Chung} (pages
%73-76), the proof follows from the method of difference equations.

%\begin{Lem} \label{3.7bb} Recalling \ref{3.7sec}, for all $\alpha$ we have
% \begin{eqnarray}
% & &
% \p^{\alpha}_x\left[T_a>T_b\right]=\frac{\sum_{i=a+1}^{x-1}\exp
%\big(S_{i}-  S_{a}\big)  +1}{\sum_{i=a+1}^{b-1}\exp\big( S_{i}-
%S_{a}  \big)+1 } \label{k1} , \\  & & \p^{\alpha}_x \left[T_a<T_b
%\right]=\frac{\sum_{i=x+1}^{b-1} \exp \big( S_{i}- S_{b} \big)
%+1}{\sum_{i=a+1}^{b-1}\exp \big( S_{i}- S_{b} \big) +1 }
%\label{k2} .
%\end{eqnarray}
%\end{Lem}

\noindent For completeness we recall an explicit expression for the mean and an upper bound for the variance
 of the local times at a certain stopping time, we can be found a proof of these elementary facts in \cite{Revesz} (page 279)
%\begin{Lem} \label{Lou} For all $\alpha$ and $i \in \Z$, we have, if $x>i$
%\begin{eqnarray}
%\Ea_i\left[
%\lo(x,T_i)\right]=\frac{\alpha_i\p_{i+1}^{\alpha}\left[T_{x}<T_{i}\right]}{
%\beta_{x}\p^{\alpha}_{x-1}\left[T_{x}>T_{i}\right]} ,
%\label{2eq65}
%\end{eqnarray}
%if $x<i$
%\begin{eqnarray}
%\Ea_i\left[ \lo(x,T_i)\right]= \frac
%{\beta_i\p_{i-1}^{\alpha}\left[T_{x}<T_{i}\right]}{
%\alpha_{x}\p^{\alpha}_{x+1}\left[T_{x}>T_{i}\right]} .
%\label{2eq68}
%\end{eqnarray}
%\end{Lem}
%\noindent Thanks to this two lemma we easily get the following one, that is used several time in this paper
\begin{Lem} \label{Lou} For all $\alpha$, 
Let $k>\tmo$ 
\begin{eqnarray}
 \Ea_{\tmo}\left[ \lo(k,T_{\tmo})\right]& =& \frac{\alpha_{\tmo}}{\beta_k}
\frac{1}{e^{S_k-S_{\tmo}}} a_{k,m_{n}} \label{A.1}, \textrm{ where} \\
a_{k,m_{n}}&=&\frac{\sum_{i=\tmo+1}^{k-1}e^{S_i}+e^{S_k}}{\sum_{i=\tmo+1}^{k-1}e^{S_i}+e^{S_{\tmo}}} \label{A.2} . \\
%b_{k,m_{n}}= \sum_{i=\tmo+1}^{k-1}e^{S_i-S_{M_{k}}}+e^{S_k-S_{M_{k}}\\
%c_{k,m_{n}}= 
 \V_{\tmo}\left[ \lo(k,T_{\tmo})\right]& \leq & 2 (\Ea_{\tmo}\left[ \lo(k,T_{\tmo})\right])^2 
\frac{e^{S_{M_{k}}-S_{m_n}}}{\beta_k} |k-m_{n}| \label{A.3}.
\end{eqnarray}
$M_{k}$ is such that $S_{M_{k}}=\max_{ m_{n}+1 \leq j \leq k-1 } S_{j}$. For $Q$-a.a. environment $\alpha$
\begin{eqnarray}
  \frac{\eta_{0}}{1-\eta_{0}}  \leq \frac{\alpha_{\tmo}}{\beta_k}  a_{k,m_{n}} \leq \frac{1}{\eta_{0}}. \label{A.4}
\end{eqnarray}
%if $k<\tmo$
%\begin{eqnarray}
% \Ea_{\tmo}\left[ \lo(k,T_{\tmo})\right]& =& \frac{\beta_{\tmo}}{\alpha_k}
%\frac{1}{e^{S_k-S_{\tmo}}}\frac{\sum_{i=k+1}^{\tmo-1}e^{S_i}+e^{S_k}}{\sum_{i=k+1}^{\tmo-1}e^{S_i}+e^{S_{\tmo}}}
%\ ,
%\end{eqnarray}
%\begin{eqnarray}
% \V_{\tmo}\left[ \lo(k,T_{\tmo})\right]& \leq &  \frac{\alpha_{\tmo}}{(\beta_k)^2}
%\frac{e^{S_{M_{k}}-S_{k}}}{e^{S_k-S_{\tmo}}} a_{k,m_{n}} b_{k,m_{n}}(1+b_{k,m_{n}}-a_{k,m_{n}})
%\end{eqnarray}
A similar result is true for $k<m_{n}$ and $\Ea_{\tmo}\left[ \lo(m_n,T_{\tmo})\right]=1$.
\end{Lem}

\noindent \textbf{Acknowledgements}
I would like to thank Nicolas Klutchnikoff for stimulating discussions.

%\noindent To show the Theorem \ref{tpslocmo}  we use the propertie
%\ref{8eq26} and the result of sub-diffusivity showed by Y. Sinai
%(Proposition \ref{lem2} page \pageref{lem2}). First let us recall
%some elementary result on Sinai's walk.
% Dans la section suivante on fait une étude locale du milieu à
%l'extérieur du voisinage $F_{p}(n)$ (point clé de la proof).
%Enfin dans le troisième paragraphe on donne la proof du result.
%\newpage
 \bibliography{thbiblio}

\vspace{1cm} \noindent
\begin{tabular}{l}
Laboratoire MAPMO - C.N.R.S. UMR 6628 - F\'ed\'eration Denis-Poisson  \\
Universit\'e d'Orl\'eans, UFR Sciences \\
B\^atiment de math\'ematiques - Route de Chartres \\
B.P. 6759 - 45067 Orl\'eans cedex 2 \\
FRANCE 
\end{tabular}

\end{document}